\documentclass{article}
\usepackage[margin=1.5in]{geometry}
\usepackage{amsmath,natbib,graphics,amsthm,amssymb,setspace,float,pdfpages,array,booktabs}
\usepackage{graphicx, epigraph}
\doublespace
\theoremstyle{definition}

\newtheorem{scenario}{Scenario}

\newcommand\numberthis{\addtocounter{equation}{1}\tag{\theequation}}

\makeatletter
\newenvironment{chapquote}[2][2em]
  {\setlength{\@tempdima}{#1}%
   \def\chapquote@author{#2}%
   \parshape 1 \@tempdima \dimexpr\textwidth-2\@tempdima\relax%
   \itshape}
  {\par\normalfont\hfill--\ \chapquote@author\hspace*{\@tempdima}\par\bigskip}
\makeatother

\title{Differences Among Noninformative Stopping Rules Are Often Relevant to Bayesian Decisions\footnote{Thanks to Samuel Fletcher, Konstantin Genin, Stephan Hartmann, Daniel Malinsky, Conor Mayo-Wilson, and the attendees at a works-in-progress session at the Munich Center for Mathematical Philosophy for comments and feedback on earlier drafts. This work was supported by a Visiting Fellowship from the Munich Center for Mathematical Philosophy.}
}
\author{Greg Gandenberger}
\date{\today}

\begin{document}
	\maketitle

\begin{chapquote}{L. J. \citet[5]{savage54}}
``It will, I hope, be demonstrated...that the superficially incompatible systems of ideas associated on the one hand with a personalistic [Bayesian] view of probability and on the other hand with the objectivistically inspired developments of the [classical] British-American School do in fact lend each other mutual support and clarification.''
\end{chapquote}

\section{Introduction}\label{sec:intro}

In classical statistics, the ``stopping rule'' that governs the decision to stop collecting data for the purpose of some inference or decision typically affects the calculations that guide that inference or decision.
In Bayesian statistics, by contrast, stopping rules typically make no difference to posterior probabilities.
For this reason, they are often taken to be irrelevant to Bayesian decisions.

However, I point out in this paper that Bayesian principles actually can require a decision-maker to attend to differences in stopping rules in a kind of situation that is common in scientific practice, in which the choice of a decision policy could influence choices about stopping rules in ways that do not align with the decision-maker's interests.
I then argue that recognizing the relevance of stopping rules to decision-making in cases of this kind is sufficient to accommodate classical statistical intuitions about stopping rules insofar as those intuitions are compelling.

The issue of the role that stopping rules should play in inference and decision has great practical significance.
It is important for decisions based on data from scientific studies to be made efficiently.
In medical research, for instance, lives may be at stake.
As a result, there has been interest in designs for clinical trials of new medical treatments that allow researchers to look at their data multiple times, each time deciding whether to declare a new treatment effective, declare it ineffective, or continue the trial.
Within a classical statistical approach, looking at the data multiple times has significant costs: it complicates trial planning and analysis, and it requires a more stringent standard for declaring the treatment to be effective on a given look.
Within a Bayesian approach, by contrast, looking at the data multiple times has no drawbacks other than the time and energy it requires: one does not need to take into account how many times one has looked or will look at an accumulating body of data in order to use that body of data to update one's beliefs in a Bayesian manner.
Bayesians often tout this feature of their approach as an advantage (e.g. \citealp{edwards+al63,berry06}), while advocates of classical methods regard it as a source of concern \citep{mayo+kruse01}.
I claim that the gap between Bayesians and frequentists on this issue is not as large as it might appear because the frequentist practice of attending to differences among stopping rules often has a Bayesian rationale.

As indicated in the epigraph, Savage had originally hoped to use Bayesian concepts to support and clarify, rather than overthrow, classical statistical ideas.
Along those lines, he claimed that the classical commitment to randomization in experimental design has a Bayesian justification, despite initial appearances to the contrary \citep[585]{savage61}.
However, he could not see a role for attention to differences among noninformative stopping rules in Bayesian procedures \citep[239]{edwards+al63}, and on the whole he came to ``[lose] faith in the devices'' of classical statistics \citep[iv]{savage72}.
Here I aim to show that Savage and other Bayesians have missed an opportunity to support and clarify, rather than overthrow, the role that stopping rules play in classical statistical inference.

\section{An Example} \label{sec:stopping}

The classical and Bayesian views about the (ir)relevance of stopping rules to inference and decision can both be made to seem unreasonable.
Consider the following pair of scenarios.
\vspace{-10 mm}
	\begin{quotation}
		\begin{scenario}
A scientist who works for a pharmaceutical company---call her ``Mammona''---is testing a new drug.
Mammona is selfish, and if the FDA approves the drug then the value of her company stock will increase.
As a result, she would prefer that the FDA approve the drug regardless of whether it is more effective than the current standard of care.
As a result, she decides to give the drug to patients and record whether or not they recover until either her funds are exhausted or the data from her trial fit the claim that the drug is more effective than the current standard of care to a particular degree.
		\end{scenario}
	\end{quotation}

\vspace{-10 mm}
	\begin{quotation}
		\begin{scenario}
In another possible world, a different pharmaceutical company scientist---call her ``Charis''---is testing the same drug as Mammona.
Charis is a humanitarian, so she would prefer that the FDA approve the drug if and only if it is more effective than the current standard of care.
As a result, she decides to give the drug to a number of patients that is fixed in advance rather than waiting for seemingly favorable results as Mammona does.
\end{scenario}
		\end{quotation}

Suppose that in these two scenarios, it just happens that Mammona and Charis treat exactly the same patients and get exactly the same outcomes.
Both scientists report their data to the FDA, as well as the ``stopping rule'' they used to decide when to stop collecting data.

It might seem obvious that the FDA should be more reluctant to approve the drug in the scenario in which it receives Mammona's results than in the scenario in which it receives Charis's results, in accordance with classical statistics.
After all, Mammona's procedure is more likely than Charis's to yield a result that fits the hypothesis that their drug is more effective than the standard of care to some degree even if that hypothesis is false.
The degree to which either of their data seem to indicate that the drug is more effective than the current standard of care will almost inevitably fluctuate somewhat as the trial proceeds, regardless of the truth.
Mammona's procedure, unlike Charis's can take advantage of this noise in the data by ``stopping on a random high.''
One might think that this fact must be taken into account in assessing the import of Mammona's reported outcomes.
If the FDA allowed procedures like Mammona's and did not account for their error characteristics in its decisions, then it would approve a greater proportion of ineffective drugs than if it insisted on fixed-sample designs like Charis's.
These considerations tend to make the Bayesian view about the relevance of stopping rules to inference and decision seem unreasonable.

On the other hand, the two scenarios involve the same hypotheses, the same patients, and the same outcomes.
How, then, could the evidential import of the data with respect to the effectiveness of the drug be different?
And if the evidential import is not different, then how could there be a difference in which inferences and decisions are warranted?
We can prefer Charis's stopping rule to Mammona's, but that is a matter of before-trial planning, rather than after-trial inference and decision.
The difference in the stopping rule in the two cases is simply a difference in the intentions of the experimenters about whether or not they would have kept going if the data had been different, which seems completely irrelevant to the import of the data their experiments actually produced.
These considerations tend to make the classical view about the relevance of stopping rules to inference and decision seem unreasonable.

How can these competing intuitions be resolved?
In the next section (Section \ref{sec:posteriors}), I explain how core Bayesian principles entail that differences among stopping rules are irrelevant to the degrees of belief one should have in the relevant hypotheses in light of the data in typical cases.
Then in Section \ref{sec:stop}, I explain why Bayesians should nevertheless take those differences into account in decision-making in circumstances that are common in scientific practice.
Finally, in Section \ref{sec:success} I argue that an appropriately nuanced Bayesian position conforms to the classical position insofar as that position is plausible and thus successfully reconciles competing intuitions about the relevance of stopping rules to inference and decision.

\section{The Irrelevance of Differences among Noninformative Stopping Rules to Posterior Probabilities}\label{sec:posteriors}

Suppose the variables $X_1,X_2,\ldots$ have a joint distribution $p(\cdot)$ that depends on the value of some parameter $\theta$, and that the values of those variables arrive in sequence.
After each observation $X_i$, you have the choice to stop observing or to continue.
That decision may be influenced by the data that have been observed so far, but it is not influenced by $\theta$ except through the data.
A stopping rule that satisfies these conditions is called \textit{noninformative}.
Given a noninformative stopping rule $s$ that gives probability $c(x^i)$ of observing another variable after observing $x^i=x_1,\ldots,x_i$, the probability density $p_s(\cdot)$ of observing exactly $x^N=x_1,\ldots,x_N$ is

\begin{align*}
p_s(x^N|\theta)&=c(\emptyset)p(x_1|\theta)c(x_1)p(x_2|x^1,\theta)\ldots c(x^{N-1})p(x_N|x^{N-1},\theta)(1-c(x^N))\\
&\propto \Pi_{i=1}^N p(x_i|x^{i-1},\theta)\\
&=p(x^N|\theta)
\end{align*}
In words, the probability that an experiment with a particular noninformative stopping rule would yield a particular body of data, given $\theta$, is proportional (as a function of $\theta$) to the probability that the variables that are actually observed would take on the values they actually take, given $\theta$.
The stopping rule plays no role in fixing that probability up to a constant of proportionality; it only affects the value of that constant.

It follows that differences among noninformative stopping rules make no difference to posterior probabilities under Bayesian conditioning \citep[37--8]{raiffa+schlaifer61}.
For Bayes's rule says that the posterior probability distribution over $\theta$ given $x^N$ is given by the following expression.

\begin{align*}
f(\theta|x^N)&=\frac{p_s(x^N|\theta)f(\theta)}{\int p_s(x^N|\theta)f(\theta)\,d\theta}\\
&=\frac{p(x^N|\theta)f(\theta)}{\int p(x^N|\theta)f(\theta)\,d\theta}
\end{align*}
where replacing $p_s$ with $p$ is possible because the constant of proportionality that the stopping rule introduces can be pulled out of the integral in the denominator and cancelled out in the numerator and the denominator.

The fact that Mammona's stopping rule is ``biased,'' stopping on evidence that favors one hypothesis over another but not vice-versa, does not make it informative in the technical sense.
Suppose her precise rule is to stop after observation $x_i$ if either $i=1000$ (at which point her funds are exhausted) or the \textit{likelihood ratio} $p(x^i|\theta_a)/p(x^i|\theta_0)$ reaches some threshold for a particular pair of values $\theta_a$ and $\theta_0$ of $\theta$.
In this case, the decision to stop after a particular observation is determined by the data that have been observed so far, without any further dependence on $\theta$.
Thus, it seems that the rule must be noninformative.

In fact, not enough has been said to determine whether Mammona's stopping rule is informative for a given observer.
It suffices to show that Mammona's stopping rule is noninformative \textit{for Mammona}, but not that it is noninformative for a third party.
Suppose the FDA regulator who receives Mammona's data---call him ``Krino''---suspects that Mammona might have run a pilot study that she did not report, and that she would not have taken the risk of using a ``biased'' stopping rule (which might, let use say, cause the FDA to disregard her data entirely) unless the pilot data had been unfavorable to the new drug.
Then Krino's probabilities for hypotheses about the effectiveness of the drug will be affected by what he learns about her stopping rule, even after he conditions on the data she reports.
\textit{For Krino}, the decision to stop at a particular point is not (epistemically) independent of $\theta$, given the data, in this scenario.
Thus, the stopping rule is informative for him.
However, it is not informative for him or anyone else simply by virtue of the fact that it is ``biased.''
Any stopping rule can be informative for someone in a scenario in which it provides him or her with information about suspected undisclosed data.
In what follows, I assume Krino knows that Mammona chooses a stopping rule by a criterion that is not sensitive to any undisclosed data she might have, so that her stopping rule is not informative for him.

Outside of special circumstances, then, Mammona and Charis's stopping rules are noninformative.
Bayesian posterior probabilities are not affected by differences among noninformative stopping rules.
This fact is often taken to imply that a Bayesian who is going to make a decision to which the value of some parameter $\theta$ is relevant need not attend to differences among stopping rules that are noninformative with respect to $\theta$.
However, this further claim does not follow and is in fact false in an important class of cases, as I explain in the next section.

\section{A Bayesian Should Attend to Differences Among Noninformative Stopping Rules in Some Cases Despite Their Irrelevance to Posterior Probabilities}\label{sec:stop}

It was shown in the previous section that differences among noninformative stopping rules are irrelevant to Bayesian posterior probabilities.
It seems to follow that they are irrelevant to Bayesian decision making.
After all, Bayesians make decisions by maximizing expected utility, and the expected utility of a decision depends only on the relevant posterior probabilities and utilities.
It seems that utility of making a particular decision (for instance, to approve a drug or not) would not typically depend on the stopping rule that was used in some experiment that has already been performed, so there does not seem to be any way for differences among noninformative stopping rules to influence Bayesian decisions.

The problem with this argument is that the utility of making a particular decision can, and often does, depend on which noninformative stopping rule is used in a relevant experiment.
To see why, consider the fact that a Bayesian typically will have preferences before an experiment is run about what stopping rule is used.
For instance, Krino might prefer Charis's stopping rule to Mammona's, or vice versa.
As in the previous section, suppose that Mammon's rule is to stop after observation $x_i$ if either $i=1000$ or $p(x^i|\theta_a)/p(x^i|\theta_0)$ reaches some threshold for a particular pair of values $\theta_a$ and $\theta_0$ of $\theta$.
Suppose that Charis's rule is to stop after observation $x_{100}$.
Suppose that the hypotheses $H_0:\theta=\theta_0$ and $H_1:\theta=\theta_a$ are known to be the only possibilities and that the best course of action from Krino's perspective is to approve Charis and Mammona's drug if $H_a$ is true but not if $H_0$ is true.
Then if Krino had to choose one of Charis and Mammona's experiments to be performed and the results reported to him, he may strictly prefer one of them over the other on grounds of expected utility.
This statement is perfectly compatible with the claim that if the two experiments were to produce the same outcome then either one would lead him to the same posterior probabilities; it arises from the fact that they give rise to different probability distributions over different sets of possible outcomes, so that, speaking informally, Krino would not necessarily expect them to be equally informative.

Now suppose that Krino engages in repeated public interactions with scientists like Mammona and Charis.
If he were to ignore differences among noninformative stopping rules, then scientists could learn this fact and judge it to be in their interests to use stopping rules like Mammona's rather than Charis's.
But perhaps Krino prefers stopping rules like Charis's.
Then it could be advantageous for him to give preferential treatment to rules like Charis's to avoid incentivizing experimental designs that he dislikes.
The crucial assumption here is that Krino is engaged in repeated interactions with scientists and that scientists are able to observe and learn from his behavior.
Under these conditions, the utility of the decision to approve a drug (or not), for instance, reflects not only the downstream consequences of administering that drug to patients (or not), but also the downstream consequences of that decision for scientists' future choices about what stopping rules to use.

I will now describe a range of scenarios in which maximizing expected utility demonstrably does require a Bayesian to attend to differences among noninformative stopping rules.
For ease of analysis, I will consider a case in which Krino is required to state up front a policy for deciding whether or not to approve Mammona's drug on the basis of her experimental results, rather than a case involving repeated interactions.
More specifically, in addition to the assumptions specified above, I assume that Krino is to specify a likelihood ratio $LR_f$ such that he will approve the drug if Mammona produces a likelihood ratio $p(x^i|\theta_a)/p(x^i|\theta_0)\geq LR_f$ in a fixed-sample experiment with some given sample size $n$, and a likelihood ratio $LR_t$ such that he will approve the drug if Mammona produces a likelihood ratio of at least $LR_t$ in a ``target-likelihood ratio'' experiment that terminates as soon as either a given sample size $m>n$ is reached or $p(x^i|\theta_a)/p(x^i|\theta_0)\geq c$ for some $c$ that Mammona specifies.
Krino's posterior probability ratio $\Pr(\theta_a|x^N)/\Pr(\theta_0|x^N)$ after conditioning on the data from either experiment depends on $x^N$ only through the likelihood ratio, so setting $LR_t>LR_f$ would amount to ``penalizing'' the target-likelihood ratio experiment relative to the fixed-sample experiment by requiring a higher posterior probability ratio to approve the drug if the first experiment is performed, and vice versa for $LR_t<LR_f$.
Mammona learns the values Krino specifies for $LR_f$ and $LR_t$, and then chooses one of the available stopping rules, administers her drug to patients in accordance with it, and reports both the stopping rule and the resulting likelihood ratio.
Krino then approves the drug or not in accordance with the rule he specified in the beginning.%
\footnote{%
The assumption that Krino binds himself to a particular course of action from the beginning of this process is needed to avoid a paradox of deterrence \citep{kavka78}.
Without this assumption, it could happen that Krino attempts to deter Mammona from performing a target-likelihood-ratio experiment by setting $LR_t$ to a value greater than the likelihood-ratio cutoff for rejecting $H_0$ that would maximize his expected utility if it were settled from the beginning that a target-likelihood-ratio experiment were to be performed.
One could then argue that if Mammona were to ``call his bluff'' by perfoming the target-likelihood-ratio experiment anway, then he should not in fact use $LR_t$  as his likelihood-ratio cutoff for rejecting $H_0$, but the cutoff that maximizes his expected utility from that point forward.
But if Mammona knows that he will proceed in this way, then his attempt to deter her from performing a target-likelihood-ratio experiment by setting $LR_t$ to a high value will be ineffective.
This thorny issue does not arise in the scenario I have specified, in which there is only one choice-point for Krino because he is bound to abide by the values for $LR_f$ and $LR_t$ that he sets in the beginning.
Krino may nevertheless update his degrees of belief in $H_0$ and $H_a$ on whatever data Mammona provides, which may influence other decisions he will make in the future.
The assumption that Mammona will perform whatever experiment maximizes her probability that Krino will approve the drug in this case rules out the possibility that she will perform that target-likelihood-ratio experiment even if he sets $LR_t=\infty$ in order to influence the decisions he might make in future trials of the same drug, for instance.
To make this assumption plausible, it may be necessary to stipulate that Mammona knows that Krino's decision will never be revisited.
}
I assume that $H_0$ and $H_a$ are simple statistical hypotheses, so that $p(x^i|\theta_a)/p(x^i|\theta_0)$ has an objective value to which the degrees of belief of an agent who satisfies the Principal Principle will conform.
I assume that Krino and Mammona's degrees of belief conform to the Principal Principle in this case, that they know this about themselves and each other, and that Mammona will choose the stopping rule that maximizes her degree of belief that Krino will approve the drug.
To avoid uninteresting technical complications, I assume that the space of possible outcomes is discrete.%
\footnote{%
The assumption that the space of possible outcomes is dicrete is always met in practice because we cannot measure or store arbitrary quantities of continuous quantities with perfect precision.
}

It is shown in Appendix \ref{app:proof} that under these circumstances, maximizing expected utility \textit{always permits} Krino to set $LR_t>LR_f$, thereby ``penalizing'' the target-likelihood ratio stopping rule.
Moreover, it \textit{requires} him to do so if and only if $\Delta_0\times P\times W>\Delta_a$, where $\Delta_0$ is the amount by which moving from a fixed-sample experiment that Krino regards as optimal to a target-likelihood ratio experiment with the same likelihood-ratio cutoff for rejection increases the probability of a likelihood ratio that exceeds that cutoff under the null hypothesis; $\Delta_a$ is the corresponding quantity under the alternative hypothesis; $P$ is the ratio of Krino's degree of belief in $H_0$ to his degree of belief in $H_a$; and $W$ is the ratio of the utility for Krino of refraining from approving the drug if $H_0$ is true to that of approving it if $H_a$ is true.
Under this condition, maximizing expected utility allows Krino to set $LR_t=\infty$, which is effectively the same as the standard practice of adopting a policy of disallowing target-likelihood-ratio experiments.
Thus, these results provide a partial Bayesian rationalization for the standard frequentist practice of simply refusing to draw any conclusions from data from experiments like Mammona's.

It is shown in Appendix \ref{app:delta} that the condition $\Delta_0\times P\times W>\Delta_a$ is satisfied exactly when $PW<(1-\beta)/[\alpha+\epsilon]$,%
\footnote{%
The conditions $\Delta_0\times P\times W>\Delta_a$ and $PW<(1-\beta)/[\alpha+\epsilon]$ might seem to be in conflict, given that the first says that $PW$ is large, while the second says that $PW$ is small.
However, $\Delta_0$ is implicitly a function of $PW$ (see Appendix \ref{app:delta}).
}
 where $\beta$ is the probability under $H_a$ that a fixed-sample experiment that Krino regards as optimal rejects $H_0$ (sometimes called that experiment's ``Type II error rate''), $\alpha$ is the corresponding probability under $H_0$ (the ``Type I error rate''), and $\epsilon$ is the amount by which a target-likelihood-ratio experiment that Krino regards as optimal falls short of the ``universal bound'' on the probability that an experiment will yield a likelihood ratio of at least $P W$ for $H_a$ against $H_0$ when $H_0$ is true.
It is always true that $PW\leq (1-\beta)/\alpha$, and as $m$ increases $\epsilon$ goes very quickly to an upper bound that approaches 0 as the amount by which the likelihood ratio is expected to overshoot $c$ if it reaches $c$ approaches zero.
Thus, $PW<(1-\beta)/[\alpha+\epsilon]$ can be expected to hold as long as $m$ is not small and the expected overshoot is not large.
It makes sense that this penalizing the target-likelihood-ratio experiment may be unnecessary if the expected overshoot is large: when it is large, the target-likelihood-ratio cutoff is less likely to be achieved, all else being equal, so running the target-likelihood-ratio experiment is less likely to be advantageous to Mammona.

The situation considered here is simpler than more realistic scenarios involving repeated interactions, but it illustrates a general point: when a decision-maker's choice of decision policy can influence choices about stopping rules in ways that do not align with his or her interests, maximizing expected utility may require him or her to attend to differences among noninformative stopping rules.
Doing so does not require any departure from standard Bayesian conditioning: the stopping rule affects the expected-utility calculations through the regulator's utilities rather than his or her probabilities, reflecting the fact that he or she has preferences about what experiments scientists will perform in the future.
There are many kinds of decision-makers who engage in repeated interactions with scientists whose interests may not align with their own, including not only government regulators, but also other scientists, journal editors, science journalists, evidence-based practitioners, and even the general public.
Any of these agents may need to attend to differences among noninformative stopping rules in their decisions in order to avoid incentivizing the use of stopping rules that they regard as undesirable.
Thus, the fact that Bayesians sometimes need to attend to differences among noninformative stopping rules in making decisions is not an idle curiosity, but a key to understanding how Bayesian methods should be used in science.

\section{Bayesians Account for Stopping Rules Insofar as They Are Relevant}\label{sec:success}

In the previous section, it was shown that a Bayesian decision-maker may need to attend to differences among noninformative stopping rules when his or her choice of decision policy could influence choices about stopping rules in ways that do not align with his or her interests.
But stopping rules play a larger role in classical statistics, where they always enter into the calculations that are used to make inferences and decisions.
In this section, I defend the Bayesian position on stopping rules against three possible objections from advocates of the classical statistical position.

First, one might claim that stopping rules are relevant to the evidential import of the data, and not just to decisions by way of utilities.
This objection is difficult to resolve because it is not clear to what the phrase ``evidential import'' refers, if anything.
A Bayesian could avoid it simply by declining to adopt any account of evidential import at all.
However, the claim that differences among noninformative stopping rules are irrelevant to evidential import is at least defensible.
That claim follows from the Likelihood Principle, which has strong axiomatic arguments in its favor \citep{birnbaum62,berger+wolpert88,gandenberger15}.
It is also supported by the informal argument given in Section \ref{sec:stopping}: a difference in stopping rules amounts to a difference in experimenters' intentions only, and it seems odd to maintain that evidential relations between data and hypotheses depend not only on the data and hypotheses themselves, but also on what the experimenters were thinking when they generated the data.
Thus, the Bayesian position on stopping rules can be defended against this first objection.

As a second objection, one might claim that stopping rules are relevant to beliefs as well as actions.
For instance, it seems that Krino should be more reluctant not only to approve the drug, but also to \textit{believe} that the drug is more effective than the current standard of care on the basis of Mammona's results rather than Charis's.

One possible response to this objection is to appeal to general arguments for Bayesian conditioning, such as those that appeal to diachronic Dutch books.
Advocates of classical statistics do not accept the use of Bayesian conditioning for typical scientific hypotheses, but they do accept it when the relevant prior probabilities are appropriately grounded in known chances or frequencies.
Thus, they take the position that differences among noninformative stopping rules are relevant for beliefs when appropriate prior probabilities are unavailable but not otherwise.
The Bayesian position that they are irrelevant to beliefs in any sense that is reducible to degrees of belief has many arguments in its favor and seems preferable to this strangely disunified view.

Another possible response to this objection is to deny that full beliefs are reducible to degrees of belief.
Under Bayesian conditioning, one's degrees of belief about a set of hypotheses are unaffected by differences among stopping rules that are noninformative with respect to those hypotheses when conditioning on data produced by those stopping rules.
However, one might have a notion of full belief that is not reducible to Bayesian degrees of belief and is spelled out, for instance, at least in part in terms of one's dispositions to act in certain ways under certain conditions.
On such an account, it might come out as true that Krino has the full belief that the drug is more effective than the current standard of care if he receives Charis's results but not if he receives Mammona's, if it turns out that maximizing expected utility would lead him to approve the drug on the basis of Charis's results but not Mammona's.

A third possible response to this objection is that our intuitions about what is relevant to belief are likely to be unreliable in this case.
Perhaps those intuitions really track relevance to decision, which generally coincides with relevance to belief but does not do so in this case.

A final possible objection is that differences among noninformative stopping rules are relevant to decisions in a broader range of circumstances than Bayesians can accommodate.
To take an extreme case,  suppose an apocalyptic event left only one survivor.
As it happens, that survivor suffers from the ailment that Charis and Mammona's drug is designed to treat.
She happens to find a supply of their drug and a report from one of their experiments.
She has to decide whether or not to take the drug.
Should she take into account the stopping rule that is described in the report?
Given that there is (let us suppose) no prospect that her decision will influence others' choices regarding stopping rules in the future, Bayesian principles seem to dictate that she should ignore the stopping rule.
And yet, one might have the intuition that she should be more hesitant to take the drug if she found Mammona's report rather than Charis's.

The Bayesian position regarding this case may seem counterintuitive, but it seems to be well supported by the standard arguments for Bayesian conditioning and the Likelihood Principle.
The most compelling arguments advocates of classical statistics give for attending to stopping rules raise the worry that if we did not attend to them, then disingenuous experimenters would be able to mislead us (see e.g
. \citealp{mayo+kruse01}).
But this worry does not arise in an apocalyptic scenario.
Moreover, it arises in precisely the cases in which Bayesian principles also yield the conclusion that differences among noninformative stopping rules are potential relevant to decisions.
Beyond a bare appeal to somewhat shaky intuitions, frequentists have no argument that we should attend more to differences among noninformative stopping rules than standard Bayesian principles require.

\section{Conclusion}

Differences among noninformative stopping rules can make a difference to a Bayesian decision-maker when his or her choice of decision policy could influence a choice regarding stopping rules in ways that do not reflect his or her interests.
This kind of situation is common in science.
As a result, Bayesian and classical statistics agree about the relevance of stopping rules to decision-making to a greater extent than is usually recognized.
Where the two positions differ, I have argued that the Bayesian position can be defended against various objections in addition to being supported by general arguments for the Likelihood Principle and for Bayesian conditioning.

\bibliography{stopping-rules}
\bibliographystyle{apalike}

\appendix

\section{Proofs About When Expected Utility Maximization Requires and When It Permits $LR_t>LR_f$}\label{app:proof}

We can distinguish four possible outcomes of Mammona's experiment for Krino.
If $H_0$ is true, then he will either incorrectly approve the drug (a ``Type I error'') or correctly decline to approve it.
If $H_a$ is true, then he will either correctly approve the drug or incorrectly decline to approve it (a ``Type II error'').
The conditional probabilities of these outcomes are fixed by Krino's choice of a likelihood-ratio cutoff for approving the drug and the experimental design.
Thus, letting $R$ indicate that $H_0$ is rejected (i.e
., the drug is approved), we can write the expected utility of a given experiment $E$ for Krino as follows.%
\footnote{%
Although the concepts of testing and rejecting hypotheses and of making Type I and Type II errors are central to frequentist but not Bayesian thinking, there is nothing un-Bayesian about representing Krino's expected utility in this way.
``Rejecting $H_0$'' here just means approving the drug, and thus does not relate directly to any notion of full belief or acceptance that cannot be reduced to Bayesian degrees of belief.
I am simply considering the utilities for Krino of two possible acts (approve the drug, do not approve the drug) in two possible states of affairs ($H_0$ is true, $H_a$ is true).
}
\begin{align*}
EU(E)&=\mbox{Pr}_E(R|H_0)\Pr(H_0)U(\mbox{Type I error})+\mbox{Pr}_E(\neg R|H_0)\Pr(H_0)U(\mbox{correct non-rejection})\\
&\qquad +\mbox{Pr}_E(R|H_a)\Pr(H_a)U(\mbox{correct rejection})+\mbox{Pr}_E(\neg R|H_a)\Pr(H_a)U(\mbox{Type II error})\\
&= \mbox{Pr}_E(R|H_a)\Pr(H_a)[U(\mbox{correct rejection}-U(\mbox{Type II error})]\\
& \qquad-\mbox{Pr}_E(R|H_0)\Pr(H_0)[U(\mbox{correct non-rejection}-U(\mbox{Type I error})]\\
& \qquad +\Pr(H_0)U(\mbox{correct non-rejection})+\Pr(H_a)U(\mbox{Type II error}) \numberthis \label{eq:eu}
\end{align*}

I assume that the prior probabilities in this equation are not affected by the experimental design.
This assumption is warranted by the stipulation that Krino knows that Mammona simply chooses the design among those available to her which maximizes her probability that Krino will approve the drug.
Thus, her choice of experimental design cannot reveal anything about her (possibly warranted) beliefs about the effectiveness of the drug.
I assume also that the utilities are not affected by the experimental design.
This assumption is warranted by the stipulation that the two simple statistical hypotheses $H_0$ and $H_a$ are the only possibilities; thus, for instance, there is no uncertainty about how effective the drug is if a Type I error is committed.
Without this assumption, a Type I error could suggest different likely sizes of departure from the null hypothesis depending on the experimental design, so that the cost of a Type I error would not be the same across experimental designs.

Fix a ``sample space'' of possible sequences of observations.
By equation \ref{eq:eu}, adding an element $x$ of that space to the ``rejection region'' of the experiment---that is, to the subset of the sample space on which $H_0$ is rejected---increases the expected utility of the experiment for Krino if and only if it increases $\mbox{Pr}_E(R|H_a)\Pr(H_a)[U(\mbox{correct rejection}-U(\mbox{Type II error})]$ for him more than it increases $\mbox{Pr}_E(R|H_0)\Pr(H_0)[U(\mbox{correct non-rejection}-U(\mbox{Type I error})]$.
It does so if and only if the likelihood ratio $\mbox{Pr}_E(x|H_a)/\mbox{Pr}_E(x|H_0)$ exceeds $P W$, where $P=\Pr(H_0)/\Pr(H_a)$ and $W=[U(\mbox{correct non-rejection}-U(\mbox{Type I error})]/[U(\mbox{correct rejection}-U(\mbox{Type II error})]$.
By assumption, $P$ and $W$ are not affected by the stopping rule.
Thus, the optimal experiment for Krino among those with a particular stopping rule uses the likelihood-ratio cutoff $P W$ for rejecting $H_0$, regardless of the stopping rule.

Now, when the sample space is discrete, as I have assumed, using $P W$ as the likelihood ratio cutoff for rejection is not unique in maximizing the expected utility for Krino of a given experiment $E$.
If an element of the sample space has exactly that likelihood ratio, then Krino is indifferent about including it in the rejection region.
In addition, there will be an open interval around $P W$ that such that varying the likelihood-ratio cutoff for rejection within that interval does not change the rejection region because there are no results in the sample space that have likelihood ratios within that interval.
This complication could be avoided through the use of continuous sample spaces, but at the cost of introducing other technical complications and the idealization of infinitely precise observations.

I have assumed that Mammona will choose the experiment that maximizes her probability that Krino approves the drug.
Thus, if she performs the target-likelihod-ratio experiment, then she will set the target likelihood ratio $c$ that suffices to end the experiment equal to the minimum value Krino specifies for the likelihood ratio $LR_t$ that would suffice for him to approve the drug if the target-likelihood-ratio experiment is performed.
Now, the target-likelihood-ratio experiment will end with a likelihood ratio for $H_a$ against $H_0$ that exceeds a particular threshold on every possible data sequence on which the fixed-sample experiment would do likewise, as well as on some possible data sequences on which the fixed-sample experiment would not do so.
Thus, if Krino sets $LR_t\leq LR_f$, then Mammona will perform the target-likelihood-ratio experiment with $c=LR_t$.%
\footnote{%
The only exception to this implication occurs when $LR_t$ and $LR_f$ are so large that the probability that the corresponding experiments would achieve them is zero.
I will ignore this scenario by restricting attention to cases in which it would not maximize Krino's expected utility to set $LR_t$ and $LR_f$ so high that they could not be achieved.
}
It follows that if the optimal fixed-sample experiment with $LR_f=PW$ has greater expected utility for Krino than the optimal target-likelihood-ratio experiment with $c=LR_t=PW$, then maximizing expected utility requires him to set $LR_f$ to some value in a neighborhood around $PW$ and $LR_t$ to some value greater than $PW$ that is sufficiently large to cause Mammona to choose the fixed-sample experiment.
Otherwise, it permits him to set $LR_f>LR_t=P W$, so that Mammona will perform the target-likelihood ratio experiment with $c=LR_t=P W$.
Thus, maximizing expected utility requires Krino to set $LR_t>LR_f$ if and only if the optimal fixed-sample experiment with $LR_f=PW$ has greater expected utility for Krino than the optimal target-likelihood-ratio experiment among those that Mammona would choose (with $c=LR_t=PW$).

Let $E_f$ be the fixed-sample experiment with $LR_f=P W$ and $E_t$ the target-likelihood-ratio experiment with $c=LR_t=P W$.
Then by the reasoning above, maximizing expected utility requires Krino to set $LR_t>LR_f$ if and only if the following condition is satisfied.

\begin{align*}
&EU(E_f)>EU(E_t)\\
&EU(E_f)-EU(E_t)>0\\
& \mbox{Pr}_{E_f}(R|H_a)\Pr(H_a)[U(\mbox{correct rejection})-U(\mbox{Type II error})]\\
&\qquad - \mbox{Pr}_{E_t}(R|H_a)\Pr(H_a)[U(\mbox{correct rejection})-U(\mbox{Type II error})]\\
&\qquad -\mbox{Pr}_{E_f}(R|H_0)\Pr(H_0)[U(\mbox{correct non-rejection})-U(\mbox{Type I error})]\\
&\qquad + \mbox{Pr}_{E_t}(R|H_0)\Pr(H_0)[U(\mbox{correct non-rejection})-U(\mbox{Type I error})] > 0\\
&[\mbox{Pr}_{E_t}(R|H_0)-\mbox{Pr}_{E_f}(R|H_0)]\Pr(H_0)[U(\mbox{correct non-rejection})-U(\mbox{Type I error})]\\
&\qquad - [\mbox{Pr}_{E_t}(R|H_a)- \mbox{Pr}_{E_f}(R|H_a)]\Pr(H_a)[U(\mbox{correct rejection}-U(\mbox{Type II error)})]> 0\\
&[\mbox{Pr}_{E_t}(R|H_0)-\mbox{Pr}_{E_f}(R|H_0)]\Pr(H_0)[U(\mbox{correct non-rejection})-U(\mbox{Type I error})]\\
&\qquad > [\mbox{Pr}_{E_t}(R|H_a)- \mbox{Pr}_{E_f}(R|H_a)]\Pr(H_a)[U(\mbox{correct rejection)}-U(\mbox{Type II error})]\\
&[\mbox{Pr}_{E_t}(R|H_0)-\mbox{Pr}_{E_f}(R|H_0)] \frac{\Pr(H_0)}{\Pr(H_a)}\frac{U(\mbox{correct non-rejection)}-U(\mbox{Type I error})}{U(\mbox{correct rejection)}-U(\mbox{Type II error})}
\\
&\qquad > [\mbox{Pr}_{E_t}(R|H_a)- \mbox{Pr}_{E_f}(R|H_a)]\\
& \Delta_0\times P\times W >\Delta_a
\end{align*}
where $\Delta_0=\mbox{Pr}_{E_t}(R|H_0)-\mbox{Pr}_{E_f}(R|H_0)$ and $\Delta_a=\mbox{Pr}_{E_t}(R|H_a)-\mbox{Pr}_{E_f}(R|H_a)$.

Therefore, maximizing expected utility requires Krino to set $LR_t>LR_f$ if and only if $\Delta_0\times P\times W >\Delta_a$.
Under this condition, it permits him to set $LR_f=P W$ and $LR_t=\infty$, because doing so would lead Mammona to perform the experiment that he regards as optimal.

Maximizing expected utility still permits Krino to set $LR_t>LR_f$ when this condition is not satisfied.
In that case, the target-likelihood ratio experiment with $c=LR_t=P W$ has at least as much expected utility for Krino as any other available experiment.
Thus, maximizing expected utility permits him to fix $LR_t$ and $LR_f$ in any way that would cause Mammona to choose that experiment.
That requires him to set $LR_t$ in a neighborhood around $P W$ and $LR_f$ sufficiently high that Mammona will choose the target-likelihood ratio experiment.
Setting $LR_f=LR_t=P W$ would accomplish this goal, but so would setting $LR_t=P W$ and $LR_f=LR_t-\epsilon$ for some $\epsilon>0$.
At the very least, because of the discreteness of the sample space, $\epsilon$ could be chosen to be sufficiently small that the fixed-sample experiment with $LR_f=P W-\epsilon$ has the same rejection region as $E_f$, so that Mammona would still choose the target-likelihood-ratio experiment.
It might also be possible to choose $\epsilon$ in a way that adds elements to the rejection region while still keeping the probability of rejecting $H_0$ for Mammona smaller on the fixed-sample experiment than on the target-likelihood-ratio experiment, but that depends on the sample space and on Mammona's probabilities for $H_0$ and $H_a$.
Therefore, maximizing expected utility always permits Krino to set $LR_t>LR_f$.

\section{Proof Regarding When $\Delta_0\times P\times W>\Delta_a$}\label{app:delta}

Let $\delta=1-\mbox{Pr}_{E_t}(R|H_a)$.
Then

\begin{align*}
\Delta_a &=\mbox{Pr}_{E_t}(R|H_a)-\mbox{Pr}_{E_f}(R|H_a)\\
&= (1-\delta)-\mbox{Pr}_{E_f}(R|H_a)\\
&= (1-\mbox{Pr}_{E_f}(R|H_a))-\delta\\
&= \mbox{Pr}_{E_f}(\neg R|H_a)-\delta\\
&= \beta-\delta
\end{align*}

Next, let $\epsilon=1/PW-\mbox{Pr}_{E_t}(R|H_0)$. I claim in the main text that $\epsilon$ is the amount by which a target-likelihood-ratio experiment that Krino regards as optimal falls short of the ``universal bound'' on the probability that an experiment will yield a likelihood ratio of at least $P W$ for $H_a$ against $H_0$ when $H_0$ is true. This claim follows from the fact that the universal bound on the probability that an experiment will yield a likelihood ratio of at least $k$ for $H_a$ against $H_0$ when $H_0$ is true is $1/k$ \citep{robbins70}.

So now we have

\begin{align*}
\Delta_0 &=\mbox{Pr}_{E_t}(R|H_0)-\mbox{Pr}_{E_f}(R|H_0)\\
&= 1/PW-\epsilon-\alpha\\
\end{align*}

and thus

\begin{align*}
\Delta_0\times P\times W &> \Delta_a\\
[1/PW-\epsilon-\alpha]PW&>\beta-\delta\\
1-PW(\alpha+\epsilon)&>\beta-\delta\\
PW(\alpha+\epsilon) &< 1-\beta+\delta\\
PW &< \frac{1-\beta+\delta}{\alpha+\epsilon}
\end{align*}

$\delta>0$, so it follows that $\Delta_0\times P\times W>\Delta_a$ if $PW < (1-\beta)/(\alpha+\epsilon)$, as stated in the main text.

It was also stated in the main text that $PW \leq (1-\beta)/\alpha$.
This claim follows from the fact that the rejection region of $E_f$ consists of the elements of the sample space on which the likelihood ratio for $H_a$ against $H_0$ is at least $PW$.
Thus, $1-\beta=\sum_i b_i$ and $\alpha=\sum_i a_i$ for some $\{a_i,b_i\}$ such that $a_i,b_i>0$ and $b_i/a_i\geq PW$ for all $i$.
As a consequence,

\begin{align*}
\frac{1-\beta}{\alpha}&=\frac{\sum_i b_i}{\sum_i a_i}\\
&=\frac{\sum_i \frac{b_i a_i}{a_i}}{\sum_i a_i}\\
&\geq \frac{\sum_i PW a_i}{\sum a_i}\\
&= PW
\end{align*}
as stated.

A result from \citet[1201]{blume08} entails that for a broad range of distributions over the sample space (namely, those that come from a full exponential family), the limit of $\epsilon$ as $m$ goes to infinity is a decreasing function of the ``expected overshoot'' and is zero when the expected overshoot is zero.
The expected overshoot is the expected size of the difference between the stopping boundary $c$ and the likelihood ratio that occurs when it reaches that boundary.
As an illustration, suppose $H_0$ says that the probability that a given patient will recover on Mammona and Charis's drug is $3/7$ (independent of any other patient outcomes), while $H_a$ says that it is $6/7$.
Then a single patient's recovery generates a likelihood ratio for $H_a$ against $H_0$ of $(6/7)/(3/7)=2$, while a single patient's failure to recover generates a likelihood ratio for $H_a$ against $H_0$ of $(1/7)/(4/7)=1/4$.
The likelihood ratio for $H_a$ against $H_0$ on the total data from multiple patients is simply the product of the individual likelihood ratios.
In this case, the expected overshoot is zero if $c$ is a positive interval power of two: if a likelihood ratio of $c=2^n$ is achieved for some $n\in \{2, 4, 8, \ldots\}$, then it is achieved exactly.
On the other hand, if $c=6$, for instance, then the expected overshoot is 2: a likelihood ratio of at least $6$ can only be achieved by passing from a likelihood ratio of $4$ to a likelihood ratio of $8$.
Cases with small expected overshoots are the cases in which the target-likelihood-ratio experiment will be particularly attractive to Mammona.
The results given here show that they are also the cases in which $\epsilon$ is small, and thus the cases in which maximizing expected utility requires Krino to penalize the target-likelihood-ratio experiment for the sake of causing Mammona not to choose it.

\end{document}